\documentclass[11pt]{amsart}
\usepackage{amssymb}
\newcommand{\cc}{\mathbb{C}}
\newcommand{\zz}{\mathbb{Z}}
\newcommand{\pp}{\mathbb{P}}
\newcommand{\dbar}{\overline{\partial}}

\DeclareMathOperator{\sign}{sgn}
\DeclareMathOperator{\supp}{supp}
\DeclareMathOperator{\opI}{{\bf I}}
\DeclareMathOperator{\opS}{{\bf S}}

\newtheorem{theorem}{Theorem}

\newtheorem{definition}[theorem]{Definition}

\title[$\dbar$-Integration in homogeneous varieties II] 
{An explicit $\dbar$-integration formula for weighted 
homogeneous varieties II,\\ forms of higher degree}

\author{J. Ruppenthal}
\author{E. S. Zeron}

\address{Fachbereich C, Bergische Universit\"at Wuppertal, 
Gau{\ss}str. 20, D-42119 Wuppertal, Germany.}
\email{ruppenthal@uni-wuppertal.de}

\address{Depto. Matem\'aticas, CINVESTAV, Apartado 
Postal 14-740, M\'exico D.F., 07000, M\'exico.}
\email{eszeron@math.cinvestav.mx}

\date{November 23, 2008}
\thanks{The second author was supported by 
Cinvestav(Mexico) and Conacyt-SNI(Mexico)}

\begin{document}
\begin{abstract}
Let $\Sigma$ be a weighted homogeneous (singular) subvariety of $\cc^n$. 
The main objective of this paper is to present a class of explicit integral 
formulae for solving the $\dbar$-equation $\omega=\dbar\lambda$ on the 
regular part of $\Sigma$, where $\omega$ is a $\dbar$-closed $(0,q)$-form 
with compact support and degree $q\geq1$. Particular cases of  these 
formulae yield $L^p$-bounded solution operators for  $1\leq p\leq \infty$ 
if $\Sigma$ is a homogeneous and pure dimensional subvariety of $\cc^n$
with an arbitrary singular locus.
\end{abstract}

\maketitle

~\\[-18mm]
\section{Introduction}

As it is well known, solving the $\dbar$-equation forms one of the 
main pillars of complex analysis, but it also has deep consequences 
on algebraic geometry, partial differential equations and other areas. 
For example, the classical Dolbeault theorem implies that the $\dbar$-equation can 
be solved in all degrees on a Stein manifold, and it is known that 
an open subset of $\cc^n$ is Stein if and only if the $\dbar$-equation can 
be solved in all degrees (on that set).
Nevertheless, it is usually not easy to produce an explicit 
operator for solving the $\dbar$-equation on a given Stein manifold, 
even if we know that it can be solved. The construction of explicit 
operators depends strongly on the geometry of the manifold on which 
the equation is considered. There exists a vast literature about this 
problem on smooth manifolds, both in books and papers (see \cite{HL,LM}, for example).

The respective Dolbeault theory on singular varieties has been developed 
only recently. Let $\Sigma$ be a singular subvariety of the space $\cc^n$
and $\omega$ a bounded $\dbar$-closed differential form on the regular part 
of $\Sigma$. Forn{\ae}ss, Gavosto and Ruppenthal have produced a general 
technique for solving the $\dbar$-equation $\omega=\dbar\lambda$ on the 
regular part of $\Sigma$, which they have successfully applied to 
varieties defined by the formula $z^m=\prod_kw_k^{b_k}$ in $\cc^n$; 
see \cite{Ga,FG} and \cite{Ru}. Acosta, Sol\'{\i}s and Zeron have 
developed an alternative technique for solving the $\dbar$-equation 
(if $\omega$ is bounded) on the regular part of any singular quotient 
variety embedded in $\cc^n$ which is generated by a finite group of unitary 
matrices, like for instance hypersurfaces in $\cc^3$ with only a Rational Double 
Point singularity; see \cite{AZ1,AZ2} and \cite{SZ}.

Nevertheless, the research on calculating explicit operators for 
solving the $\dbar$-equation $\omega=\dbar\lambda$ on the regular part 
of singular subvarieties $\Sigma\subset\cc^n$ is still at a very early 
state; the techniques mentioned in the previous paragraph do not produce 
useful explicit formulae. Ruppenthal and Zeron have proposed explicit operators 
for calculating solutions $\lambda$ if $\Sigma$ is a weighted 
homogeneous variety and $\omega$ is a $\dbar$-closed $(0,1)$-differential 
form with compact support; see \cite{RZ}. The weighted homogeneous 
varieties are analysed, for they are a main model for classifying the 
singular subvarieties of $\cc^n$. A detailed analysis of the weighted 
homogeneous varieties is done in Chapter~2--\S4 and Appendix B of 
\cite{Di}. The main objective of the present paper is to improve the explicit 
operators originally developed in \cite{RZ} for calculating solutions 
$\lambda$ to the $\dbar$-equation $\omega=\dbar\lambda$ on the regular 
part of any weighted homogeneous variety $\Sigma$ if $\omega$ 
is a $\dbar$-closed $(0,q)$-differential form with compact support and 
degree $q\geq1$.
Furthermore, we produce $\dbar$-solution operators 
with $L^p$-estimates for $1\leq{p}\leq\infty$ if $\Sigma$ is homogeneous
with an arbitrary singular locus.

\begin{definition}\label{whs}
Let $\beta\in\zz^n$ be a fixed integer vector with strictly positive 
entries $\beta_k\geq1$. A holomorphic polynomial $Q(z)$ on $\cc^n$ is said 
to be \textbf{weighted homogeneous} of degree $d\geq1$ with respect to 
$\beta$ if the following equality holds for all $s\in\cc$ and $z\in\cc^n$:
\begin{eqnarray}
\label{homo}&&Q(s^{\beta}*z)\,=\,
s^d\,Q(z),\quad\hbox{with~the~action:}\\
\label{action}&&s^{\beta}*(z_1,z_2,...,z_n)\,:=\,
(s^{\beta_1}z_1,s^{\beta_2}z_2,...,s^{\beta_n}z_n).
\end{eqnarray}

An algebraic subvariety $\Sigma$ in $\cc^n$ is said to be 
\textbf{weighted homogeneous} with respect to $\beta$ if $\Sigma$ 
is the zero locus of a finite number of weighted homogeneous polynomials 
$Q_k(z)$ of (maybe different) degrees $d_k\geq1$, but all of them with 
respect to the same fixed vector $\beta$.
\end{definition}

\newpage
Let $\Sigma\subset\cc^n$ be any subvariety. We use the following notation 
along this paper. The regular part $\Sigma^*=\Sigma_{reg}$ is the complex 
manifold consisting of the regular points of $\Sigma$, and it is always 
endowed with the induced metric, so that $\Sigma^*$ is a Hermitian 
submanifold in $\cc^n$ with corresponding volume element $dV_\Sigma$ and 
induced norm $|\cdot|_\Sigma$ on the Grassmannian $\Lambda{T}^*\Sigma^*$. 
Thus, any Borel-measurable $(0,q)$-form $\omega$ on $\Sigma^*$ admits a 
representation $\omega=\sum_Jf_Jd\overline{z_J}$, where the coefficients 
$f_J$ are Borel-measurable functions on $\Sigma^*$ which satisfy the 
inequality $|f_J(z)|\leq |\omega(z)|_\Sigma$ for all points $z\in\Sigma^*$ 
and multi-indexes $|J|=q$. Notice that such a representation is by no 
means unique. We refer to Lemma 2.2.1 in \cite{Ru} for a more detailed 
treatment of that point. For $1\leq p<\infty$, we also introduce the 
$L^p$-norm of a measurable $(0,q)$-form $\omega$ on an open set 
$U\subset\Sigma^*$ via the formula:
\begin{eqnarray*}
\|\omega\|_{L^p_{0,q}(U)} &:=&\bigg(\int_{U} 
|\omega|^p_\Sigma\,dV_\Sigma\bigg)^{1/p}.
\end{eqnarray*}

We can now present the main result of this paper. We assume that the 
$\dbar$-differentials are calculated in the sense of distributions, 
for we work with Borel-measurable functions.

\begin{theorem}[\textbf{Main}]\label{main}
Let $\Sigma$ be a weighted homogeneous subvariety of $\cc^n$ with 
respect to a given vector $\beta\in\zz^n$, where $n\geq2$ and all 
entries $\beta_k\geq1$. Consider the class of all $(0,q)$-forms 
$\omega$ given by $\sum_Jf_Jd\overline{z_J}$, where $q\geq1$, the 
coefficients $f_J$ are all Borel-measu\-rable functions in $\Sigma$, 
and $z_1,..., z_n$ are the Cartesian coordinates of $\cc^n$. Let 
$\sigma\geq-q$ be any fixed integer. The operator  $\opS^\sigma_q$ below 
is well defined on $\Sigma$ for all forms $\omega$ which are essentially 
bounded and have compact support,
\begin{eqnarray}\label{fng-1}
&&\opS^\sigma_q\omega(z):=\sum_{|J|=q}\frac{\aleph_J}{2\pi{i}}
\int_{u\in\cc}f_J(u^{\beta}*z)\,\frac{u^\sigma(\overline{u^{\beta_J}})
d\overline{u}\wedge{du}}{\overline{u}\,(u-1)}\\
\label{fng-2}&&\hbox{with}\quad\aleph_J=
\sum_{j\in{J},\,K=J\setminus\{j\}}
\frac{\beta_j\overline{z_j}\,d\overline{z_K}}{\sign(j,K)}
\quad\hbox{and}\quad\beta_J=\sum_{j\in{J}}\beta_j.
\end{eqnarray}

Notice that the multi-indexes $J$ and $K$ are both ordered in an 
ascending way and that $\sign(j,K)$ is the sign of the permutation 
used for ordering the elements of the $q$-tuple $(j,K)$ into an 
ascending way. Finally, the form $\opS^\sigma_q(\omega)$ is a solution 
of the $\dbar$-equation $\omega=\dbar\opS^\sigma_q(\omega)$ on the regular 
part of $\Sigma\setminus\{0\}$, whenever $\omega$ is also $\dbar$-closed 
on the regular part of $\Sigma\setminus\{0\}$. 
\end{theorem}

The origin of $\cc^n$ is in general a singular point of $\Sigma$ 
according to Definition~\ref{whs}, so that the regular parts of $\Sigma$ 
and $\Sigma\setminus\{0\}$ coincide. We will prove Theorem~\ref{main} 
in Section~\ref{section:main} of this paper. Similar techniques and a slight 
modification of equations (\ref{fng-1}) and (\ref{fng-2}) can also be used 
for producing a $\dbar$-solution operator with $L^p$-estimates on 
homogeneous subvarieties with arbitrary singular locus.

\begin{theorem}[\textbf{$L^p$-Estimates}]\label{thm:l2} 
Let $\Sigma$ be a pure $d$-dimensional homogeneous (cone) subvariety of 
$\cc^n$, where $n\geq2$ and each entry $\beta_k=1$ in Definition~\ref{whs}.
Fix a real number $1\leq{p}\leq\infty$ and an integer $1\leq{q}\leq{d}$. 
Consider the class $L^p_{0,q}(\Sigma)$ of all $(0,q)$-forms $\omega$ 
given by $\sum_Jf_Jd\overline{z_J}$, where the coefficients $f_J$ are 
all $L^p$-integrable functions in $\Sigma$, and $z_1,...,z_n$ are the 
Cartesian coordinates of $\cc^n$. Choose $\sigma\in\zz$ to be the 
smallest integer such that
\begin{equation}\label{eq:l1}
\sigma\,\geq\,\frac{2d-2}p+1-q.
\end{equation}

The operator $\opS^\sigma_q(\omega)$ below is well defined 
almost everywhere on $\Sigma$ for all forms $\omega$ which lie 
in $L^p_{0,q}(\Sigma)$ and have compact support on $\Sigma$:
\begin{eqnarray}\label{eq:l2}
&&\opS^\sigma_q \omega (z):=\sum_{|J|=q}\frac{\aleph_J}
{2\pi{i}}\int_{u\in\cc}f_J(uz)\frac{u^\sigma\,\overline{u^q}\,
d\overline{u}\wedge{du}}{\overline{u}\,(u-1)},\\
\nonumber&&\hbox{where}\quad\aleph_J=\sum_{j\in{J},\,K=J\setminus\{j\}}
\frac{q\,\overline{z_j}\,d\overline{z_K}}{\sign(j,K)}.
\end{eqnarray}

The form $\opS^\sigma_q(\omega)$ is a solution of the $\dbar$-equation 
$\omega=\dbar\opS^\sigma_q(\omega)$ on the regular part of 
$\Sigma\setminus\{0\}$, whenever $\omega$ is also $\dbar$-closed on the 
regular part of $\Sigma\setminus\{0\}$. Finally, assuming that the support 
of $\omega$ is contained in an open ball $B_R$ of radius $R>0$ and centre at 
the origin, there exists a strictly positive constant $C_\Sigma(R,\sigma)$ 
which does not depend on $\omega$ and such that:
\begin{equation}\label{eq:l3}
\big\|\opS^\sigma_q(\omega)\big\|_{L^p_{0,q-1}(\Sigma\cap{B_R})}
\leq{}C_\Sigma(R,\sigma)\cdot\|\omega\|_{L^p_{0,q}(\Sigma)}.
\end{equation} 
\end{theorem}

The case $p=\infty$ in the previous theorem is a corollary 
of Theorem~\ref{main} because the formulae \eqref{eq:l2} and 
\eqref{fng-1} coincide in the homogeneous case (where all coefficients 
$\beta_J=q$). We will give the full proof of Theorem~\ref{thm:l2} in 
Section~\ref{section:L2} of the present paper.

\newpage
The obstructions to solving the $\dbar$-equation with $L^p$-estimates
on subvarieties of $\cc^n$ are not completely understood in general.
An $L^2$-solution operator (for forms with non-compact support) is
only known in the case where $\Sigma$ is a complete
intersection\footnote{More precisely: a Cohen-Macaulay space.}
of pure dimension $\geq 3$ with only isolated singularities. This
operator was constructed by Forn{\ae}ss, {\O}vrelid and Vassiliadou in
\cite{FOV2} via an extension theorem for $\dbar$-cohomology groups
originally presented by Scheja \cite{Sch}. Usually, the $L^p$-results
come with some obstructions to the solvability of the $\dbar$-equation.
Different situations have been analysed in the works of Diederich,
Forn{\ae}ss, {\O}vrelid, Ruppenthal and Vassiliadou: It is shown that
the $\dbar$-equation is solvable with $L^p$-estimates for forms lying
in a closed subspace of finite codimension of the vector space of all
the $\dbar$-closed $L^p$-forms if the variety has only isolated
singularities \cite{DFV,Fo,FOV2,OV,Ru1}. Besides, in the paper \cite{FOV1},
the $\dbar$-equation is solved locally with some weighted $L^2$-estimates
for forms which vanish to a sufficiently high order on the (arbitrary)
singular locus of the given varieties.

There is a second line of research about the $\dbar$-operator on complex
projective varieties (see \cite{PaSt1, PaSt2} for the state of the art
and further references). Though that area has clearly a lot in common with
the topic of $\dbar$-equations on analytic subvarieties of $\cc^n$, it
is a somewhat different theory because of the strong global tools (like
Serre duality) which cannot be used in the (local) situation of Stein
spaces (due to the lack of compactness).

Since the estimates in Theorem~\ref{thm:l2} are given only for homogeneous 
varieties, we finally propose in Section~\ref{final} of this paper a useful 
technique for generalising the estimates in Theorem~\ref{thm:l2}, so as to
consider weighted homogeneous subvarieties instead of homogeneous ones.

\newpage
\section{Proof of Main Theorem}\label{section:main}

We need the following result. The notation $L^1_{p,q}$ stands for the
class of all the $(p,q)$-forms with $L^1$-integrable coefficients, so
that the differentials are calculated in the sense of distributions.

\begin{theorem}\label{first-solution}
Let $U\subset\cc^m$ be open, $2\leq{q}\leq{m}$, and
$\omega\in{L}^1_{0,q}(U)$ be a $\dbar$-closed form with compact 
support along the first coordinate $z_1$, that is, such that
$\supp(\omega)\cap{F_y}$ is compact in $U\cap{F_y}$ for all fibres
$F_y=\cc{\times}\{y\}$ with $y\in\cc^{m-1}$. Assume that $\omega$ 
is given by:
$$\omega\,=\sum_{|J|=q,\,1\notin{J}}[a_J]d\overline{z_J}\,+
\sum_{|K|=q-1,\,1\notin{K}}[a_{1,K}]d\overline{z_1}\wedge
d\overline{z_K},$$ 
where the multi-indexes $J$ and $K$ are both ordered 
in an ascending way. The following operator
\begin{eqnarray*}
\opS_q(\omega)&:=&\sum_{|K|=q-1,\,1\notin{K}}
\opI[a_{1,K}]d\overline{z_K}\,,\quad\hbox{with}\\
\opI{f}(z_1,...,z_m)&:=&\frac{1}{2\pi{i}}\int_{t\in\cc}
f(t,z_2,...,z_m)\frac{d\overline{t}\wedge{dt}}{t-z_1},
\end{eqnarray*}
is defined almost everywhere in $U$ and 
satisfies $\omega=\dbar\opS_q(\omega)$.  
\end{theorem}

Notice that $\opS_q(\omega)$ is well defined in $U$ if $\omega$ is
essentially bounded and has compact support along the first coordinate
$z_1$. 

\begin{proof} It is clear that the restrictions $(a_{1,K})|_{F_y}$ are
all $L^1$-integrable on the intersections $U\cap{F_y}$, for almost every
fibre $F_y$, so that $\eta:=\opS_q(\omega)$ is defined almost everywhere
in $U$; see Appendix~B of \cite{Ra} or \cite{LM,Ru}. We only need to
show that $\dbar\eta=\omega$. The assumption $\dbar\omega=0$ implies
that the following equation holds for every multi-index $|J|=q$ with
$1\notin{J}$,
\begin{equation}\label{closed}
\frac{\partial[a_J]}{\partial\overline{z_1}}\,=
\sum_{j\in{J},\,K=J\setminus\{j\}}\sign(j,K)
\frac{\partial[a_{1,K}]}{\partial\overline{z_j}}.
\end{equation}

The function $\sign(j,K)$ is the sign of the permutation used for
ordering the elements of the $q$-tuple $(j,K)$ into an ascending way. 
A direct application of the inhomogeneous Cauchy-Integral Formula in one
complex variable and the fact that $\omega$ has compact support along
the first coordinate yield the following identity for every multi-index
$|K|=q{-}1$ with $1\notin{K}$:
$$\dbar(\opI[a_{1,K}])=[a_{1,K}]d\overline{z_1}\wedge{d}\overline{z_K}
+\sum_{j\notin{K},\,j\neq1}\opI\bigg[\frac{\partial[a_{1,K}]}
{\partial\overline{z_j}}\bigg]d\overline{z_j}\wedge{d}\overline{z_K},$$
and so we have that:
\begin{eqnarray*}
&&\partial\opS_q(w)=\sum_{|K|=q-1,\,1\notin{K}}[a_{1,K}]d\overline{z_1}
\wedge{d}\overline{z_K}\,+\sum_{|J|=q,\,1\notin{J}}\opI(b_J)d\overline{z_J}\\
&&\quad\hbox{with}\quad{b}_J\,:=\sum_{j\in{J},\,K=J\setminus\{j\}}
\sign(j,J)\frac{\partial[a_{1,K}]}{\partial\overline{z_j}}.
\end{eqnarray*}

Recall that the multi-indexes $J$ and $K$ are both ordered in an
ascending way and $\sign(j,K)$ is the sign of the permutation used for
ordering the elements of the $q$-tuple $(j,K)$ into an ascending way.
Equation~(\ref{closed}) implies that $\partial\opS_q(w)$ is equal to
$\omega$, because $a_K$ has compact support along the first coordinate,
and so:
$$\opI(b_J)\,=\,\opI\bigg(\frac{\partial[a_J]}
{\partial\overline{z_1}}\bigg)\,=\,a_J.$$
\end{proof}

We may now proceed with the proof of the main theorem.

\begin{proof}\textbf{[Main Theorem~\ref{main}]}. We follow the proof
originally presented in~\cite{RZ}, so that we only point out the main
points. Let $\{Q_k\}$ be the set of polynomials which define the
algebraic variety $\Sigma$ as its zero locus. The definition of weighted
homogeneous varieties implies that the polynomials $Q_k(z)$ are all
weighted homogeneous with respect to the same fixed vector $\beta$.
Equation~(\ref{homo}) automatically yields that every point
$s^{\beta}*z$ lies in $\Sigma$ for all $s\in\cc$ and $z\in\Sigma$, 
and so each coefficient $f_J(\cdot)$ in equation~(\ref{fng-1}) is well
evaluated in $\Sigma$. Moreover, the coefficients $\beta_k\geq1$ and
$\beta_J\geq{q}$, for all index $k$ and multi-index $J$ of degree $q$.
Fixing any point $z\in\Sigma$, the given hypotheses imply that the
following Borel-measu\-rable functions are all essentially bounded 
and have compact support in $\cc$,
$$u\,\mapsto\,f_J(u^{\beta}*z).$$

Hence, the operator $\opS^\sigma_q(\omega)$ in
(\ref{fng-1})--(\ref{fng-2}) is well defined on $\Sigma$ for 
each fixed integer $\sigma\geq-q$ and all forms $\omega$ which 
are essentially bounded and have compact support. We shall prove 
that $\opS^\sigma_q(\omega)$ is also a solution of the equation
$\omega=\dbar\opS^\sigma_q(\omega)$ if the $(0,q)$-form $\omega$ is
$\dbar$-closed. We may suppose, without loss of generality and because
of the given hypotheses, that the regular part of $\Sigma$ does not
contain the origin. Let $\xi\neq0$ be any fixed point in the regular
part of $\Sigma$. We may suppose by simplicity that the first entry
$\xi_1\neq0$, and so we define the following mapping and subvariety:
\begin{equation}
\label{eqn1}\begin{array}{rcl}
\eta(y)&:=&(y_1/\xi_1)^\beta*(\xi_1,y_2,y_3,...,y_n),
\quad\hbox{for}\quad{y}\in\cc^n,\\
Y&:=&\{\widehat{y}\in\cc^{n-1}:Q_k(\xi_1,\widehat{y})=0\,\forall\,k\}.
\end{array}\end{equation}

The action $s^{\beta}*z$ was given in (\ref{action}). We have that
$\eta(\xi)=\xi$, and that the following identities hold for all
$s\in\cc$ and $\widehat{y}\in\cc^{n-1}$, recall equation~(\ref{homo})
and the fact that $\Sigma$ is the zero locus of the polynomials
$\{Q_k\}$:
\begin{equation}
\label{eqn2}\begin{array}{rcl}
Q_k(\eta(s,\widehat{y}))&=&(s/\xi_1)^{d_k}
\,Q_k(\xi_1,\widehat{y})\quad\hbox{and}\\
\eta(\cc^*\times{Y})&=&\{z\in\Sigma:z_1\neq0\}.
\end{array}\end{equation}

The symbol $\cc^*$ stands for $\cc\setminus\{0\}$. The mapping $\eta(y)$
is locally a biholomorphism whenever the first entry $y_1\neq0$. Whence,
the point $\xi$ lies in the regular part of the variety $\cc\times{Y}$,
because $\xi=\eta(\xi)$ also lies in the regular part of $\Sigma$ and
$\xi_1\neq0$. Thus, we can find a biholomorphism
$$\pi=(\pi_2,...,\pi_n):U\rightarrow{Y}\subset\cc^{n-1}$$
defined from an open and bounded domain $U$ in $\cc^m$ onto an open 
set in the regular part of $Y$, such that $\pi(\zeta)$ is equal to
$(\xi_2,...,\xi_n)$ for some $\zeta\in{U}$. Consider the following
holomorphic mapping defined for all points $s\in\cc$ and $x\in{U}$,
\begin{equation}\label{eqn3} 
\Pi(s,x)\,:=\,s^{\beta}*(\xi_1,\pi(x))
\,=\,\eta(s\xi_1,\pi(x))\,\in\,\Sigma. 
\end{equation}

The image $\Pi(\cc\times U)$ will be known as a {\bf generalised cone}
from now on. Notice that $\Pi(\cc^*\times{U})$ lies in the regular 
part of $\Sigma$, for $\pi(U)$ is contained in the regular part of $Y$. 
The mapping $\Pi(s,x)$ is locally a biholomorphism whenever $s\neq0$,
because $\eta$ is also a local biholomorphism for $y_1\neq0$; and the
image $\Pi(1,\zeta)$ is equal to $\xi$. Hence, recalling the form
$\omega$ and the operator $\opS^\sigma_q(\omega)$ defined in
(\ref{fng-1})--(\ref{fng-2}), we only need to prove that the 
pull-back $\Pi^*\omega$ is equal to the differential
$\dbar\Pi^*\opS^\sigma_q(\omega)$ inside $\cc^*\times{U}$, in order to
conclude that the $\dbar$-equation $\omega=\dbar\opS^\sigma_q(\omega)$
holds in a neighbourhood of $\xi$ in $\Sigma$. We can use equations
(\ref{action}) and (\ref{eqn3}) in order to calculate the pull-back
$\Pi^*\omega$ when $\omega$ is given by $\sum_Jf_Jd\overline{z_J}$. 
To simplify the notation, let $\pi_1(x):=\xi_1$ for all $x\in U$, 
so that $d\pi_1=0$.
\begin{eqnarray}\label{eqn4}
&&[\Pi^*\omega](s,x)\,=\,\sum_{|J|=q}f_J(\Pi(s,x))\,
\overline{s^{\beta_J}}\bigwedge_{j\in{J}}d\overline{\pi_j(x)}\,+\\
\nonumber&&+\sum_{|J|=q,\,j\in{J}}\frac{f_J(\Pi)\beta_j
\overline{s^{\beta_J-1}\pi_j(x)}}{\sign(j,J\setminus\{j\})}\,d\overline{s}
\wedge\!\bigwedge_{k\in{J}\setminus\{j\}}d\overline{\pi_k(x)}.
\end{eqnarray}

Recall that $\beta_J=\sum_{j\in{J}}\beta_j\geq{q}\geq-\sigma$, 
the multi-index $J$ is ordered in an ascending way, and
$\sign(\alpha_1,...,\alpha_q)$ is the sign of the permutation used for
ordering the elements of the $q$-tuple $(\alpha_1,...,\alpha_q)$ into an
ascending way. The given hypotheses on $\omega$ yield that the pull-back
$\Pi^*\omega$ is $\dbar$-closed and bounded in $\cc^*\times{U}$, and so
it is also $\dbar$-closed in $\cc\times{U}$; see Lemma~4.3.2 in
\cite{Ru} or Lemma~(2.2) in \cite{SZ}. The same argument applies 
to the $\dbar$-closed and essentially bounded form
\begin{eqnarray}\label{eq:s-1}
s^{\sigma+1}\,[\Pi^*\omega](s,x)\,\in\,L^\infty_{0,q}(\cc\times{U}).
\end{eqnarray}

The open set $U$ is bounded in $\cc^m$. Thus, it follows from \eqref{eqn4} 
and \eqref{eq:s-1}, by the use of Lemma~7.2.2 in \cite[p.~186]{Ru} or 
Lemma~3.6 in \cite{Ru2}, that $s^{\sigma}\Pi^*\omega$ is both 
$L^1_{0,q}(\cc\times U)$ and $\dbar$-closed in $\cc\times{U}$. It is 
easy to see that each coefficient $f_J(\Pi(s,x))$ has compact support 
with respect to the first coordinate $s$, so that we can apply 
Theorem~\ref{first-solution} to $t^\sigma[\Pi^*\omega](t,x)$ and 
calculate the form:
\begin{eqnarray}\label{eqn5}
&&\opS_q(t^\sigma\Pi^*\omega)=\sum_{|J|=q}\frac{\Theta_J}{2\pi{i}}
\int_{t\in\cc}f_J(\Pi(t,x))\frac{t^{\sigma}(\overline{t^{\beta_J}})
d\overline{t}\wedge{dt}}{\overline{t}\,(t-s)}\\
\nonumber&&\hbox{with}\quad\Theta_J=\sum_{j\in{J}}
\frac{\beta_j\overline{\pi_j(x)}}{\sign(j,J\setminus\{j\})}
\bigwedge_{k\in{J}\setminus\{j\}}d\overline{\pi_k(x)}.
\end{eqnarray}

Theorem~\ref{first-solution} implies that 
$$s^\sigma\,[\Pi^*\omega](s,x)=\dbar\opS_q(t^\sigma\,[\Pi^*\omega](t,x)).$$

Hence, we only need to verify that the form
$\opS_q(t^\sigma\Pi^*\omega)/s^\sigma$ is equal to the pull-back
$\Pi^*\opS^\sigma_q(\omega)$ of the form defined in~(\ref{fng-1}), in
order to conclude that $\omega=\dbar\opS^\sigma_q(\omega)$, as desired.
We begin by calculating the pull-back $\Pi^*\aleph$ of the differential
form $\aleph_J$ given in~(\ref{fng-2}). Notice that
$\pi_1(x)\equiv\xi_1$, so that $d\pi_1=0$, and recall 
equations (\ref{action}) and (\ref{eqn3}).
\begin{eqnarray}\label{eqn6}
&&\quad\Pi^*\aleph_J\,=\,\sum_{j\in{J}}\frac{\beta_j
\overline{s^{\beta_J}\pi_j(x)}}{\sign(j,J\setminus\{j\})}
\bigwedge_{k\in{J}\setminus\{j\}}d\overline{\pi_k(x)}\,+\\
\nonumber&&+\sum_{j,k\in{J},\,j\neq{k}}\frac{\beta_j\beta_k
\overline{s^{\beta_J-1}\pi_j(x)\pi_k(x)}\quad{d}\overline{s}}
{\sign(j,J\setminus\{j\})\sign(k,J\setminus\{j,k\})}\wedge
\bigwedge_{i\in{J}\setminus\{j,k\}}d\overline{\pi_i(x)}
\end{eqnarray}

Suppose that $J=(\alpha_1,...,\alpha_a,j,\dot{\alpha}_1,...,
\dot{\alpha}_b,k,\ddot{\alpha}_1,...,\ddot{\alpha}_c)$, then:
\begin{eqnarray*}
\sign(j,J\setminus\{j\})\sign(k,J\setminus\{j,k\})&=&(-1)^a(-1)^{a+b},\\
\sign(k,J\setminus\{k\})\sign(j,J\setminus\{j,k\})&=&(-1)^{a+b+1}(-1)^a.
\end{eqnarray*}

So that the last sum in equation~(\ref{eqn6}) vanishes, 
and so the pull-back $\Pi^*\aleph_J$ is identically equal 
to $\overline{s^{\beta_J}}\Theta_J$, with $\Theta_J$ defined
in~(\ref{eqn5}). Finally, we can calculate the pull-back of 
the form $\opS^\sigma_q(\omega)$ given in~(\ref{fng-1}), 
notice that $\Pi(us,x)$ is equal to $u^\beta*\Pi(s,x)$,
\begin{equation}\label{eqn7}
\Pi^*\opS^\sigma_q(w)=\sum_{|J|=q}\frac{\overline{s^{\beta_J}}\Theta_J}
{2\pi{i}}\int_{u\in\cc}f_J(\Pi(us,x))\frac{u^\sigma(\overline{u^{\beta_J}})
d\overline{u}\wedge{du}}{\overline{u}\,(u-1)}.
\end{equation}

The change of variables $t=us$ yields that the form
$\opS_q(t^\sigma\Pi^*\omega)/s^\sigma$ in~(\ref{eqn5}) is equal to the
identity above, and so $\omega=\dbar\opS^\sigma_q(\omega)$, as desired.
\end{proof}

\section{$L^p$-Estimates}\label{section:L2}

We prove Theorem~\ref{thm:l2} in this section. Recall that $\Sigma$ is 
a pure $d$-dimensional homogeneous (cone) subvariety of $\cc^n$ with 
arbitrary singular locus, so that $n\geq2$ and each entry $\beta_k=1$ 
in Definition~\ref{whs}. Moreover, given a fixed real number 
$1\leq{p}\leq\infty$ and an integer $1\leq{q}\leq{d}$, we consider 
the class $L^p_{0,q}$ of all $(0,q)$-forms $\omega$ expressed as 
$\sum_Jf_Jd\overline{z_J}$, where the coefficients $f_J$ are all 
$L^p$-integrable functions in $\Sigma$, and $z_1,...,z_n$ are the 
Cartesian coordinates of $\cc^n$. Assume that the support of each form 
$\omega\in{L}^p_{0,q}$ is contained in the open ball $B_R$ of radius 
$R>0$ and centre at the origin. Fix $\sigma\in\zz$ to be the smallest 
integer such that
\begin{equation}\label{eqn8}
\sigma\,\geq\,\frac{2d-2}p+1-q.
\end{equation}

We begin by showing that $\opS^\sigma_q$ 
in~ (\ref{eq:l2}) defines a bounded operator
\begin{eqnarray}\label{eqn9}
\opS^\sigma_q\,:\,L^p_{0,q}(\Sigma\cap B_R)
\,\rightarrow\,L^p_{0,q-1}(\Sigma\cap B_R),
\end{eqnarray}
where:
\begin{eqnarray}\label{eq:l4}
&&\opS^\sigma_q\omega (z)=\sum_{|J|=q}\frac{\aleph_J}
{2\pi{i}}\int_{u\in\cc}f_J(uz)\frac{u^\sigma\,\overline{u^q}
\,d\overline{u}\wedge{du}}{\overline{u}\,(u-1)}\\
\nonumber&&\hbox{and}\quad\aleph_J=\sum_{j\in{J},\,K=J\setminus\{j\}}
\frac{q\,\overline{z_j}\,d\overline{z_K}}{\sign(j,K)}.
\end{eqnarray}

Recall that the multi-indexes $J$ and $K$ are both ordered in an ascending 
way and that $\sign(j,K)$ is the sign of the permutation used for ordering 
the elements of the $q$-tuple $(j,K)$ into an ascending way. Notice that 
the case $p=\infty$ in (\ref{eqn9}) is a corollary of Theorem~\ref{main},
because the formulae (\ref{fng-1}) and (\ref{eq:l4}) coincide when the 
variety $\Sigma$ is homogeneous (so that all $\beta_J=q$). Hence, we can
suppose from now on that $p<\infty$, and we only need to prove that the 
following inequality holds for every multi-indexes $|J|=q$ and $j\in{J}$ 
in order to conclude that (\ref{eqn9}) and (\ref{eq:l3}) holds,
\begin{equation}\label{eqn10}
\int_{z\in\Sigma\cap{B_R}}\bigg|\overline{z_j}\!\int_{|u|<R/\|z\|}
\!\frac{f_J(uz)u^{\sigma+q}}{(u-1)\,u}dV_\cc\bigg|^pdV_\Sigma
\lesssim\|\omega\|^p_{L^p_{0,q}(\Sigma)}.
\end{equation}

Notice that the support of $f_J(z)$ is contained in the open ball $B_R$ 
of radius $R$, and that $dV_{\cc}$ and $dV_{\Sigma}$ are the respective 
volume forms on $\cc$ and $\Sigma$. Further, we may use the variable $u$
instead of its complex conjugate $\overline{u}$ because we work under an 
absolute value sign. Let $\delta<1$ be any fixed real number. It is easy 
to deduce the existence of a finite real constant $M_1$ such that the 
following inequalities hold for all complex numbers $\hat{t}$ and 
$\hat{w}$:
$$\int_{|w|<R}\frac{dV_{\cc}(w)}{|\hat{t}|^\delta\,|w-\hat{t}|}
\leq{M_1}\quad\hbox{and}\quad\int_{|t|<R}\frac{dV_{\cc}(t)}
{|t|^\delta\,|\hat{w}-t|}\leq{M_1}.$$

Hence, the generalised Young inequality for convolution integrals 
yields that the modified Cauchy-Pompeiu formula defines an 
$L^p$-bounded operator; see for example the Appendix~B of \cite{Ra} 
with $s=1$ and $\delta<1$:
\begin{equation}\label{eq:cif1}
\int_{|t|<R}\bigg|\int_{|w|<R}\!\frac{h(w)dw\wedge{d}\overline{w}}
{(w-t)\,|t|^\delta}\bigg|^pdV_\cc(t)\lesssim\int_{|t|<R}|h(t)|^pdV_\cc.
\end{equation}

Moreover, let $\widetilde{\Sigma}$ be the projective variety 
associated to $\Sigma$ in the space $\cc\pp^{n-1}$, for $\Sigma$ is a 
pure $d$-dimensional homogeneous subvariety of $\cc^n$. We also use the 
fact that any integral on $\Sigma$ can be decompose as a pair of nested 
integrals on $\cc$ and $\widetilde{\Sigma}$, that is:
$$\int_{z\in\Sigma}\Phi(z)\,dV_\Sigma(z)=
\int_{[z]\in\widetilde{\Sigma}}\int_{t\in\cc}\Phi(\dot{z}t)
|t|^{2d-2}\,dV_\cc(t)dV_{\widetilde{\Sigma}}([z]);$$
where on the right hand side, $\dot{z}\in\Sigma$ is any representative 
of $[z]\in\widetilde{\Sigma}$ with $\|\dot{z}\|=1$. Finally, since 
$\sigma\in\zz$ is the smallest integer which satisfies (\ref{eqn8}),
we have that the following constant
\begin{equation}\label{eqn11}
\delta:=\sigma+q-1+\frac{2-2d}p\quad
\hbox{satisfies}\quad0\leq\delta<1.
\end{equation}

We can now use the results presented in the paragraphs above 
in order to calculate (\ref{eqn10}) and (\ref{eq:l3}), with 
the change of variables $w=ut$,
\begin{eqnarray*}
&&\int_{z\in\Sigma\cap{B_R}}\bigg|\overline{z_j}\int_{|u|<R/\|z\|}
\!\frac{f_J(uz)u^{\sigma+q}}{(u-1)\,u}dV_\cc\bigg|^pdV_\Sigma\\
&&\leq\int_{[z]\in\widetilde{\Sigma}}\int_{|t|<R}|t|^p\bigg|
\int_{|u|<R/|t|}\!\frac{f_J(u\dot{z}t)u^{\sigma+q}}{(u-1)\,u}
dV_\cc\bigg|^p|t|^{2d-2}dV_{\cc}dV_{\widetilde{\Sigma}}\\
&&=\int_{[z]\in\widetilde{\Sigma}}\int_{|t|<R}\bigg|\int_{|w|<R}
\frac{f_J(w\dot{z})w^{\sigma+q-1}}{(w-t)t^{\sigma+q-2}}\cdot\frac
{dV_\cc}{|t|^2}\bigg|^p|t|^{2d-2+p}dV_{\cc}dV_{\widetilde{\Sigma}}\\
&&=\int_{[z]\in\widetilde{\Sigma}}\int_{|t|<R}\bigg|\int_{|w|<R}
\!\frac{f_J(w\dot{z})w^{\sigma+q-1}}{(w-t)\,|t|^\delta}
dV_\cc\bigg|^pdV_{\cc}dV_{\widetilde{\Sigma}}\\
&&\lesssim\int_{[z]\in\widetilde{\Sigma}}\int_{|t|<R}\big|f_J
(t\dot{z})t^{\sigma+q-1}\big|^pdV_{\cc}dV_{\widetilde{\Sigma}}\\
&&=\int_{[z]\in\widetilde{\Sigma}}\int_{|t|<R}\big|f_J(t\dot{z})
\big|^p|t|^{p\delta+2d-2}dV_{\cc}dV_{\widetilde{\Sigma}}\\
&&\leq\int_{z\in\Sigma\cap{B_R}}\!|f_J(z)|^pR^{p\delta}
dV_{\Sigma}\,\lesssim\|f_J\|^p_{L^p(\Sigma)}\leq
\|\omega\|_{L^p_{0,q}(\Sigma)}^p<\infty.
\end{eqnarray*}

We have used (\ref{eqn11}) and (\ref{eq:cif1}) with 
$h(w)=f_J(w\dot{z})w^{\sigma+q-1}$. That completes 
the proof of equations \eqref{eqn10} and \eqref{eq:l3}.

Finally, notice that the operators $\opS^\sigma_q(\omega)$ given in 
(\ref{fng-1}), (\ref{eq:l2}) and (\ref{eq:l4}) are all the same, because 
the coefficients $\beta_J=q$ for every multi-index $|J|=q$. Therefore, 
we can show that the operator $\opS^\sigma_q(\omega)$ satisfies the 
differential equation $\omega=\dbar\opS^\sigma_q(\omega)$ following step 
by step the proof presented in Section~\ref{section:main}. We only need 
to rewrite the pull-back given in (\ref{eqn4}), which is $\dbar$-closed 
in the product $\cc^*\times{U}$,
\begin{eqnarray}\label{eqn44}
&&[\Pi^*\omega](u,x)\,=\,\sum_{|J|=q}f_J(\Pi(u,x))\,
\overline{u^q}\bigwedge_{j\in{J}}d\overline{\pi_j(x)}\,+\\
\nonumber&&+\sum_{|J|=q,\,j\in{J}}\frac{f_J(\Pi)\,q\,
\overline{u^{q-1}\pi_j(x)}}{\sign(j,J\setminus\{j\})}\,d\overline{u}
\,\wedge\bigwedge_{k\in{J}\setminus\{j\}}d\overline{\pi_k(x)}.
\end{eqnarray}

And we must show that $u^\sigma[\Pi^*\omega](u,x)$ lies in 
$L^1_{0,q}(\cc\times{U})$, where $U$ is a bounded domain in $\cc^m$. 
Thus, we have that $u^\sigma\Pi^*\omega$ is also $\dbar$-closed in 
$\cc\times{U}$, because of Lemma~7.2.2 in \cite[p.~186]{Ru} or 
Lemma~3.6 in \cite{Ru2}. We can then apply Theorem~\ref{first-solution} 
and follow step by step the proof presented in Section~\ref{section:main} 
from equation~(\ref{eqn5}) to the end of that section.

Recall that the integer $\sigma\geq\frac{2d-2}p+1-q$. We begin showing
that the form $u^\sigma\Pi^*\omega$ lies in $L^p_{0,q}(\cc\times{U})$.
Notice that $\Pi(u,x)$ is equal to $u(\xi_1,\pi(x))$ because each entry 
$\beta_k=1$ in~(\ref{action}) and~(\ref{eqn3}). It is easy to calculate 
the pull-back of the volume form $dV_\Sigma$:
\begin{eqnarray}\label{eqn17}
\Pi^*dV_\Sigma=\sum_{|I|=|J|=d}\beta_{I,J}(z)dz_I
\wedge{d}\overline{z_J}\Big|_{z=u(\xi_1,\pi(x))}&&\\
\nonumber=\Theta(x)|u|^{2d-2}\big[du\wedge
d\overline{u}\big]\wedge\bigwedge_{k=1}^{d-1}
\big[dx_k\wedge{d}\overline{x_k}\big].&&
\end{eqnarray}

Recall that $x$ lies in the bounded open set $U\subset\cc^{d-1}$. 
Since $\Sigma$ is a pure $d$-dimensional homogeneous (cone) subvariety 
of $\cc^n$, the coefficients $\beta_{I,J}(z)$ are all invariant under 
the transformations $z\mapsto{uz}$, and so $\Theta(x)$ only depends on 
the vales of $\pi(x)$ and all its partial derivatives (it is constant 
with respect to $u$). The fact that $\Pi$ is a biholomorphism from 
$\cc^*\times{U}$ onto its image also implies that $\Theta$ cannot vanish. 
Hence, choosing a smaller set $U$ if it is necessary, we can suppose that 
$|\Theta|$ is bounded from below by a constant $M_2>0$. 

On the other hand, since $\Pi(u,x)=u(\xi_1,\pi)$ and the support of 
each $f_J(z)$ is contained in a ball of radius $R>0$ and centre at the 
origin, we have that every $f_J(\Pi(u,x))$ vanishes if $|u\xi_1|>R$. Thus, 
equation (\ref{eqn44}) and the analysis done in the paragraphs above imply 
that the form $u^\sigma\Pi^*\omega$ lies in $L^p_{0,q}(\cc\times{U})$, 
because the following inequalities hold for every multi-index $J$ and 
exponent $b=0,1$:
\begin{eqnarray*}
\int_{\cc\times{U}}|u^{\sigma+q-b}f_J(\Pi)|^pdV_{\cc\times{U}}&\lesssim
&\int_{\cc\times{U}}|f_J(\Pi)|^p\Theta(x)|u|^{2d-2}dV_{\cc\times{U}}\\
=\;\int_{\Pi(\cc\times U)}|f_J|^pdV_\Sigma&\leq&
\|\lambda\|^p_{L^2_{0,1}(\Sigma)}\,<\,\infty.
\end{eqnarray*}

Recall that $p(\sigma+q-b)\geq2d-2$, because of the hypothesis 
(\ref{eq:l1})--(\ref{eqn8}) imposed to $\sigma\in\zz$. Finally, the 
support of $\Pi^*\omega$ is bounded in $\cc\times{U}$, because $U$ is 
bounded and each $f_J(\Pi(u,x))$ vanishes if $|u\xi_1|>R$. Thus, we have 
that the form $u^\sigma\Pi^*\omega$ is $L^1_{0,q}$ and $\dbar$-closed 
in $\cc\times{U}$; see for example Lemma~7.2.2 in \cite[p.~186]{Ru} or 
Lemma~3.6 in \cite{Ru2}. We can then apply Theorem~\ref{first-solution}  
and follow step by step the proof presented in Section~\ref{section:main} 
from equation~(\ref{eqn5}) to the end of that section, in order to 
conclude that the operator $\opS^\sigma_q(\omega)$ satisfies the 
differential equation $\omega=\dbar\opS^\sigma_q(\omega)$, as we want.

\section{Weighted Homogeneous Estimates}\label{final}

We want to close this paper presenting a useful technique for generalising 
the estimates given in Theorem~\ref{thm:l2}, so as to consider weighted 
homogeneous subvarieties instead of cones. Let $\Sigma\subset\cc^n$ be a 
weighted homogeneous subvariety defined as the zero locus of a finite set 
of polynomials $\{Q_k\}$. Thus, the polynomials $Q_k(z)$ are all weighted 
homogeneous with respect to the same vector $\beta\in\zz^n$, and each 
entry $\beta_k\geq1$. Define the following holomorphic mapping:
\begin{equation}\label{eqn20}
\Phi:\cc^n\to\cc^n,\quad\hbox{with}\quad\Phi(x)
=(x_1^{\beta_1},x_2^{\beta_2},...,x_n^{\beta_n}).
\end{equation}

It is easy to see that each polynomial $Q_k(\Phi)$ is homogeneous, 
and so the subvariety $X\subset\cc^n$ defined as the zero locus of 
$\{Q_k(\Phi)\}$ is a cone. Consider a $(0,q)$-form $\omega$ given by 
the sum $\sum_Jf_Jd\overline{z_J}$, where the coefficients $f_J$ are 
all Borel-measu\-rable functions with compact support in $\Sigma$, 
and $z_1,..., z_n$ are the Cartesian coordinates of $\cc^n$. We 
may follow two different paths in order to solve the equation 
$\overline{\partial}\lambda=\omega$.  We may calculate the pull-back:
$$\Phi^*\omega\,=\,\sum_{|J|=q}f_J(\Phi(x))\bigg[\prod_{j\in{J}}
\beta_j\,\overline{x_j^{\beta_j-1}}\bigg]\,d\overline{x_J};$$
and then apply Theorems~\ref{main} and~\ref{thm:l2} on the 
cone $X$, so as to get the following operators:
\begin{eqnarray}\label{eqn21}
&&\opS^\sigma_q(\Phi^*\omega):=\sum_{|J|=q}\frac{\widehat{\aleph_J}}
{2\pi{i}}\int_{u\in\cc}f_J(\Phi(ux))\frac{u^\sigma(\overline
{u^{\beta_J}})d\overline{u}\wedge{du}}{\overline{u}\,(u-1)}\\
\nonumber&&
\hbox{with}\quad\widehat{\aleph_J}=\sum_{j\in{J},\,K=J\setminus\{j\}}
\frac{\beta_j\,\overline{x_j^{\beta_j}}}{\sign(j,K)}\bigg[\prod_{k\in{K}}
\beta_k\overline{x_k^{\beta_k-1}}\bigg]\,d\overline{x_K}.
\end{eqnarray}

On the other hand, we may use the main Theorem~\ref{main} 
on the weighted homogeneous variety $\Sigma$, so as to get:
\begin{eqnarray}\label{eqn22}
&&\opS^\sigma_q(\omega):=\sum_{|J|=q}\frac{\aleph_J}{2\pi{i}}
\int_{u\in\cc}f_J(u^{\beta}*z)\frac{u^\sigma(\overline{u^{\beta_J}})
d\overline{u}\wedge{du}}{\overline{u}\,(u-1)}\quad\hbox{with}\\
\nonumber&&\aleph_J=\sum_{j\in{J},\,K=J\setminus\{j\}}
\frac{\beta_j\overline{z_j}\,d\overline{z_K}}{\sign(j,K)}
\quad\hbox{and}\quad\beta_J=\sum_{j\in{J}}\beta_j.
\end{eqnarray}

We can easily verify that $\Phi^*\opS^\sigma_q(\omega)$ is equal to 
$\opS^\sigma_q(\Phi^*\omega)$. The main problem is that $\Phi^*\omega$ 
may not lie necessarily in $L^p_{0,q}(X)$ for $p<\infty$.

\bibliographystyle{plain}

\end{document}